\documentclass[a4paper,11pt]{article}
\usepackage{color,xcolor,ucs}
\usepackage[top=0.3in, bottom=0.5in, left = 0.65in, right = 0.65in]{geometry}
\usepackage[linkcolor=black,colorlinks=true,urlcolor=blue]{hyperref}
\usepackage{mathtools}

\usepackage{color,xcolor,ucs}
\usepackage{mathtools}   \usepackage{tikz} 
\usepackage{ amssymb }
\usepackage{extarrows} 
\usepackage{pgf,tikz}
\usepackage{float}
\usetikzlibrary{positioning}
\usetikzlibrary{shapes.geometric}
\usetikzlibrary{shapes.misc}
\usetikzlibrary{arrows}
\usepackage{caption}
\usepackage{mathrsfs}
\usetikzlibrary{arrows,shapes,automata,backgrounds,petri,positioning}
\usetikzlibrary{decorations.pathmorphing}
\usetikzlibrary{decorations.shapes}
\usetikzlibrary{decorations.text}
\usetikzlibrary{decorations.fractals}
\usetikzlibrary{decorations.footprints}
\usetikzlibrary{shadows}
\usetikzlibrary{calc}
\usetikzlibrary{spy}
\usepackage{amsmath}
\usepackage{array}
\usepackage{ amssymb }
\usepackage{braket}
\usepackage{qcircuit}
\usepackage{soul}
\usepackage{braket} 
\usepackage{relsize}

\usepackage{amsmath}
\usepackage{ amssymb }
\usepackage{braket}
\usepackage{qcircuit}
\usepackage{soul}
\usepackage{braket}

\title{{\LARGE The phase transition for the Gaussian free field is sharp }}
\author{Pete Rigas}
\date{}

\begin{document}

\maketitle

\begin{abstract}
  We prove that the phase transition for the Gaussian free field (GFF) is sharp. In comparison to a previous argument due to Rodriguez in 2017 which characterized a $0-1$ law for the Massive Gaussian Free Field by analyzing crossing probabilities below a threshold $h_{**}$, we implement a strategy due to Duminil-Copin and Manolescu in 2016, which establishes that two parameters are equal, one of which encapsulates the probability of obtaining an infinite connected component under free boundary conditions, while the other encapsulates the natural logarithm of the probability of obtaining a connected component from the origin to the box of length $n$ which is also taken under free boundary conditions. We quantify the probability of obtaining  crossings in easy and hard directions, without imposing conditions that the graph is invariant with respect to reflections, in addition to making use of a differential inequality adapted for the GFF. The sharpness of the phase transition is characterized by the fact that below a certain height parameter of the GFF, the probability of obtaining an infinite cluster a.s. decays exponentially fast, while above the parameter, the probability of obtaining an infinite cluster occurs a.s. with good probability.  \footnote{\textbf{MSC Class}: 60K35; 82D02}
\end{abstract}

\section{Introduction}

\subsection{Overview}

\noindent The Gaussian Free Field (GFF) is a mathematical object that continues to attract great attention from mathematicians and physicists alike. On the mathematical front, several recent works have established connections with percolation, whether it be existence of a phase transition {\color{blue}[3]}, delocalization of the height function for the six-vertex model under sufficiently flat boundary conditions {\color{blue}[2]}, adaptations of the argument for sloped boundary conditions in the six-vertex model, with applications to the Ashkin-Teller, generalized random-cluster, and $(q_{\sigma}, q_{\tau})$-cubic models {\color{blue}[11]}, construction of an IIC-type limit {\color{blue}[12]}, and, more generally, analysis of crossing probabilities in several models {\color{blue}[5,7,10]}. To contribute to rapid developments in the field, we implement a strategy due to Duminil-Copin and Manolescu in {\color{blue}[3]}, which the authors leverage for demonstrating the sharpness of the phase transition for the random-cluster model, which can be used for studying other models which satisfy similar properties.

\subsection{Statements of previous results for the random-cluster model}

We provide an overview of results for establishing the sharpness of the phase transition for the random-cluster model, as provided in {\color{blue}[4]}, and then describe similar properties which are satisfied by the GFF. Given a finite graph $G \equiv \big( V_G, E_G \big)$, for edge weight $p \in [0,1]$, and cluster-weight $q >0$, the random-cluster \textit{probability measure} of sampling a \textit{random-cluster configuration} $\omega \in \big\{ 0, 1 \big\}^{E_G}$, under boundary conditions $\xi$, is defined by,

\begin{align*}
  \phi^{\xi}_{p,q,G} \big( \omega \big) \equiv \phi \big( \omega \big) = \frac{p^{o( \omega )} \big( 1 - p \big)^{c ( \omega )} q^{k ( \omega )}}{Z}  \text{ } \text{ , } 
\end{align*}

\noindent where $o( \omega )$ denotes the number of open edges, $c( \omega)$ denotes the number of closed edges, $k ( \omega )$ denotes the number of clusters, and $Z \equiv Z \big( p , q , G\big)$ denotes the partition function which is a normalizing constant so that $\phi$ is a probability measure. The boundary conditions of the \textit{random-cluster} measure are understood as a partition of the vertices. In the following statements below, abbreviate $\phi^0_{p,q} \big( \cdot \big) \equiv \phi^0 \big( \cdot \big)$. To study the connectivity properties between two points $x$ and $y$ of the graph, which we denote with,

\begin{align*}
   \big\{    x \longleftrightarrow y       \big\} \text{ } \text{ , }
\end{align*}

\noindent equipped with $\phi$, for $q \geq 1$ on a planar, locally-finite doubly periodic connected graph $\mathscr{G}$ that is invariant under reflections with respect to the line $\big\{ (0,y) , y \in \textbf{R} \big\}$, \textbf{Theorem} \textit{1.1} of {\color{blue}[4]} asserts the existence of some $p_c \equiv p_c \big(\mathscr{G}\big)$ for which:

\begin{itemize}
    \item[$\bullet$] Given $p < p_c$, there exists $c \equiv c \big( p , \mathscr{G}\big)$ such that a path of open edges exists between $x,y \in \mathscr{G}$, in which, \begin{align*}
    \phi \big( x \longleftrightarrow y \big) \leq \mathrm{exp} \big(  - c \big| x- y \big| \big)  \text{ } \text{ . } 
\end{align*}

\item[$\bullet$] Given $p > p_c$, there exists a.s. an infinite open cluster under $\phi \big( \cdot \big)$, in which, \begin{align*}
    \phi \big( | C \big( x,y \big)| = + \infty \big) > 0     \text{ } \text{ , } 
\end{align*}

\noindent where $C \big( x,y \big)$ denotes the cluster between $x,y \in \mathscr{G}$.
\end{itemize}

\noindent To demonstrate that such a sharp phase transition exists for this $p_c$, additional properties of $\phi$ are used, including,

\begin{itemize}
    \item[$\bullet$] \textit{FKG inequality}: Given two increasing events $A,B$, and boundary conditions $\xi$, \begin{align*}
    \phi^{\xi} \big( A \cap B \big) \geq \phi^{\xi} \big( A \big) \phi^{\xi} \big( B \big)     \text{ } \text{ , } 
\end{align*}

\item[$\bullet$] \textit{Domain Markov Property (DMP)}: Given boundary conditions $\xi$, one has the equality, \begin{align*}
  \phi \big( \omega|_{G} = \cdot \big| \chi \big( \omega\big)    \equiv \xi \big) \equiv \phi^{\xi} \big( \cdot \big)    \text{ } \text{ , } 
\end{align*}

\noindent for a \textit{random-cluster configuration} $\chi \big( \omega \big)$.

\item[$\bullet$] \textit{Comparison between boundary conditions (CBC)}: Given two pairs of boundary conditions $\xi_1$ and $\xi_2$, with $\xi_1 \leq \xi_2$, $q \geq 1$, and $p$, for an increasing event $A$, $\phi^{\xi_1} \big( A \big)  \leq \phi^{\xi_2} \big( A  \big)$.

\item[$\bullet$] \textit{Comparison between edge parameters}: Given two edge parameters $p_1$ and $p_2$, with $p_1 \leq p_2$, for boundary conditions $\xi$ and $q \geq 1$, and an increasing event $A$, $\phi_{p_1, q} \big( A \big) \equiv \phi_{p_1} \big( A \big) \leq \phi_{p_2} \big( A \big) \equiv \phi_{p_2 , q} \big( A \big)$.
\end{itemize}

\noindent With the FKG, SMP, CBC, and comparison between edge parameters properties, explicitly the threshold $p_c$ for which the statement of \textbf{Theorem} \textit{1.1} holds is given by the probability of obtaining an infinite \textit{open} path, with,

\begin{align*}
  p_c \equiv \mathrm{inf} \big\{    p \in \big( 0 ,1 \big) :        \phi^0 \big( x \longleftrightarrow + \infty \big) > 0         \big\}   \text{ } \text{ , } 
\end{align*}

\noindent while another closely related threshold is explicitly given by the probability of obtaining an open path to the boundary of the box of size length $n$, $\partial \Lambda_n$, with,

\begin{align*}
 \widetilde{p_c} \equiv   \mathrm{sup} \big\{    p \in \big( 0 ,1 \big) :  \underset{n \longrightarrow + \infty}{\mathrm{lim}}  - \frac{1}{n} \mathrm{log} \big[     \phi^0        \big( 0 \longleftrightarrow \partial \Lambda_n \big)         \big]    \big\}  \text{ } \text{ . } 
\end{align*}

\noindent To exhibit that $p_c$ and $\widetilde{p_c}$ are equal, in {\color{blue}[4]} the authors employ three steps, in which crossing probabilities in hard and easy directions are quantified. 

\bigskip

\noindent In addition to all of the aforementioned quantities, \textit{differential inequalities} for the random-cluster model play a role, the first of which states, for the same increasing event $A$ and boundary conditions, that,

\begin{align*}
     \frac{\mathrm{d}}{\mathrm{d} p} \phi^{\xi}_{p,q} \big(      A \big) \geq c \text{ }    \phi^{\xi}_{p,q,q} \big(      A \big) \big( 1 -  \phi^{\xi}_{p,q} \big(      A \big) \big)  \text{ } \mathrm{log} \big(    \frac{m_{A,p}^{-1}}{2}              \big)   \text{ } \text{ , } 
\end{align*}

\noindent for some strictly positive $c$, where,

\begin{align*}
 m_{A,p} \equiv    \underset{e \in E_G }{\mathrm{max}}  \big(  \phi^{\xi}_{p,q} \big( A       \big| w (e) = 1  \big)   -  \phi^{\xi}_{p,q} \big(  A       \big| w(e) = 0  \big)   \big)                  \text{ } \text{ . }
\end{align*}

\noindent The second differential inequality states, for $H_A$ the \textit{hamming distance} between $\omega$ and $A$, that,

\begin{align*}
     \frac{\mathrm{d}}{\mathrm{d} p} \mathrm{log} \big( \phi^{\xi}_{p,q} \big( A  \big) \big) \geq \frac{\phi^{\xi}_{p,q} \big( H_A  \big)}{p \big( 1 - p \big) }      \text{ } \text{ . } 
\end{align*}

\noindent In the next section, we state analogues to each property, if they exists, for the GFF.

\subsection{GFF properties}

\noindent For the case of the GFF, the discrete and continuous version of the field satisfies the FKG inequality (both {\color{blue}[9]} and {\color{blue}[13]} contain statements of FKG for the GFF and models closely related to the GFF). To define the GFF, one must specify a mean and covariance function. First, the mean of the GFF is taken to be zero, while the covariance function is of the form, 

\begin{align*}
 \textbf{E} \big[ \phi_u \phi_v] = G \big( u , v \big) =  \textbf{E}_u \big[ \int_{0}^{\infty} \textbf{1}_{{S_t}=v} \text{ } \mathrm{d} t ]  \text{ } \text{ , } 
\end{align*}

\noindent for fields $\phi_u$ and $\phi_v$ respectively centered at $u$ and $v$, where $\textbf{E} \big[ \cdot \big]$ denotes the expectation with respect to the GFF law $\textbf{P}_h \big[ \cdot \big] \equiv \textbf{P}\big[ \cdot \big]$ (an expression for the law is also provided in {\color{blue}[13]}). For GFF level set percolation, under $\textbf{P} \big[ \cdot \big]$, the connectivity event that two points $x$ and $y$ in the graph are connected above a height threshold $h$ is,

\begin{align*}
  \big\{  x    \overset{\geq h}{\longleftrightarrow}     y  \big\}  \text{ } \text{ . } 
\end{align*}

\noindent In addition to the FKG property, the GFF satisfies the following properties:

\begin{itemize}
    \item[$\bullet$] \textit{GFF Strong Markov Property (SMP)} ({\color{blue}[1]}, \textbf{Theorem} \textit{8}): For any random connected compact connected subset $K \subsetneq G$, conditionally upon the filtration $\mathcal{F}_K$, one has the following equality,

    \begin{align*}
   \big\{    \phi_v : v \in G \backslash K       \big\}       \overset{\mathrm{d}}{\equiv} \big\{    \textbf{E} \big[ \phi_v | \mathcal{F}_K \big] + \phi_v : v  \in G \backslash K         \big\} \text{ } \text{ . } 
    \end{align*}

    \item[$\bullet$] \textit{Comparison between height parameters of the free field}. For two height parameters $h_1 \leq h_2$, and any $x,y \in G$, $\textbf{P}\big[ x \overset{\geq h_2}{\longleftrightarrow} y \big] \leq \textbf{P} \big[ x \overset{\geq h_1}{\longleftrightarrow} y \big]$.
\end{itemize}

\noindent Equipped with the FKG, Strong Markov, and comparison between height parameters properties, about the height threshold $h \equiv 0$, the sharpness of the phase transition for the GFF can be captured through the following two regimes of behavior.

\bigskip

\noindent \textbf{Theorem} \textit{1} (\textit{sharpness of the phase transition for the GFF}). For $G = \big( V, E \big)$, with $x, y \in G$, one has two possible behaviors:

\begin{itemize}
    \item[$\bullet$] Given $h < 0$, there exists $c \equiv c \big( h , G \big)$ such that,

    \begin{align*}
        \textbf{P} \big[  x \overset{\geq h}{\longleftrightarrow}  y \big]   \leq \mathrm{exp} \big( - c  \big| x - y \big| \big)    \text{ } \text{ . } 
    \end{align*}

    \item[$\bullet$] Given $h > 0$, there exists a.s. an infinite open cluster under $\textbf{P} \big[ \cdot \big]$, in which,

     \begin{align*}
        \textbf{P} \big[  \big|          C \big( x , y \big)        \big| = + \infty  \big]      > 0           \text{ } \text{ , } 
    \end{align*}

    \noindent where $C\big( x , y \big)$ denotes the cluster between $x,y \in G$.
\end{itemize}

\noindent We introduce the parameters,

\begin{align*}
    h_c \equiv \mathrm{inf} \big\{ h \in \big( -\infty , 0 \big) : \textbf{P} \big( x \overset{\geq h}{\longleftrightarrow} + \infty \big) > 0   \big\}      \text{ } \text{ , } 
\end{align*}

\noindent and, for $\Lambda_n \subsetneq G$,

\begin{align*}
     \widetilde{h_c} \equiv \mathrm{sup} \big\{   h \in \big( - \infty , 0 \big) :        \underset{n \longrightarrow + \infty}{\mathrm{lim}}   - \frac{1}{n} \mathrm{log} \big[  \textbf{P} \big(  0 \overset{\geq h}{\longleftrightarrow} \partial \Lambda_n    \big)       \big]         \big\}   \text{ } \text{ . } 
\end{align*}

\noindent For an increasing event $A$, the GFF satisfies a differential inequality, which takes the form,

 \begin{align*}
     \frac{\mathrm{d}}{\mathrm{d} h}   \textbf{P} \big[     A                   \big]  \geq  c^{\prime} \textbf{P} \big[ A \big] \big( 1  - \textbf{P} \big[ A \big] \big) \mathrm{log} \big( \frac{\mathcal{I}^{-1}_{A,h}}{2} \big)    \text{ } \text{ , } 
    \end{align*}

\noindent for some strictly positive $c^{\prime}$, where the influence term in the logarithm is,

\begin{align*}
  \mathcal{I}_{A,h} \equiv \textbf{P} \big( A \big|   \phi_x \geq h      \big) - \textbf{P} \big( A \big|    \phi_x < h     \big)  \text{ } \text{ . }
\end{align*}

\subsection{Paper organization}

\noindent In the remaining sections of the paper, we exhibit that the two height parameters defined in the previous section are equal. This exhibits the sharpness of the free field, as crossing events occurring in the hard direction are shown to correspond to the connectivite probabilities decaying exponentially fast. Above the critical height threshold, the remaining possibility is shown to hold in the final section.

\section{Crossing probabilities in the easy direction}

\noindent Introduce,

\begin{align*}
\underset{n \longrightarrow + \infty }{\mathrm{lim\text{ }  inf}}   \text{ }   \textbf{P} \big( \mathscr{V}   \mathscr{C} \big( n , 2 n \big)    \big)          \text{ } \text{ , } 
\end{align*}

\noindent for the \textit{ vertical crossing event} $\mathscr{V}\mathscr{C} \big( n , 2 n \big)$ of height $\geq h$. In the statement below, the fact that,

\begin{align*}
\underset{n \longrightarrow + \infty }{\mathrm{lim\text{ }  inf}}   \text{ }   \textbf{P} \big( \mathscr{V}   \mathscr{C} \big( n , 2 n \big)    \big)        \longrightarrow 0   \text{ } \text{ , } 
\end{align*}

\noindent implies exponential decay of a connectivity event to the boundary of a finite volume of length $n$.

\bigskip

\noindent \textbf{Proposition} \textit{1} (\textit{limit infimum of vertical crossings from n to 2n implies exponential decay}). Fix some $h_c > 0$. If $h < h_c$, there exists an infinite volume measure for which,

\begin{align*}
\underset{n \longrightarrow + \infty }{\mathrm{lim\text{ }  inf}}   \text{ }   \textbf{P} \big( \mathscr{V}   \mathscr{C} \big( n , 2 n \big)    \big)        \longrightarrow 0   \text{ } \text{ , } 
\end{align*}

\noindent then there exists some $c \equiv c \big( h \big)$ so that,

\begin{align*}
     \textbf{P} \big( 0 \overset{\geq h}{\longleftrightarrow }\partial \Lambda_n     \big) \leq \mathrm{exp} \big( - c \big| x- y \big| \big)    \text{ } \text{ . } 
\end{align*}

\noindent To show that the infinite volume measure in proposition above exists,  introduce the following lemma.

\bigskip

\noindent \textbf{Lemma} \textit{1} (\textit{exponential decay in the infinite volume measure}). For the same $h_c$ as in \textbf{Proposition} \textit{1}, there exists strictly positive $\kappa$, and $h < h_c$, for which the infinite volume measure satisfies,

\begin{align*}
  \textbf{P}_h \big[ 0 \overset{\geq h}{\longleftrightarrow}   \partial \Lambda_n            \big]  \text{ } \text{ . } 
\end{align*}

\noindent There exists $c \equiv c \big( h \big)$ such that, for any $n \geq 0$,

\begin{align*}
    \textbf{P}_h \big[ 0 \overset{\geq h}{\longleftrightarrow}   \partial \Lambda_n            \big] \leq \mathrm{exp} \big( - c \big| x -  y \big| \big)     \text{ } \text{ . } 
\end{align*}

\noindent For two height parameters, the inequality below relates how the probability of $\mathscr{V}\mathscr{C}\big( n , 2n \big)$ occurring differs.

\bigskip

\noindent \textbf{Lemma} \textit{2} (\textit{the probability of a vertical crossing from n to 2n occurs is }). Fix $h_1 \geq h_2$. For any $N \geq n$,

\begin{align*}
       \textbf{P}_{h_2} \big[   \mathscr{V} \mathscr{C} \big( n , 2n \big)    \big]      \leq \mathrm{exp} \bigg(           - \big( h_1 - h_2 \big)   \frac{N}{n} \big(           1 - \textbf{P}_{h_2} \big[ \mathscr{V}\mathscr{C} \big( n , 2n \big) ]                    \big)^{2 \frac{N}{n}}                            \bigg) \text{ } \text{ . } 
\end{align*}

\noindent \textit{Proof of Lemma 2}. To demonstrate that an exponential upper bound of the form given above holds, observe that in order for $\mathscr{V}\mathscr{C} \big( n , 2 n \big)$ to occur, either,

\begin{align*}
  \textbf{P}_{h_2} \bigg[  [    kn    , \big( k +2 \big) n ]  \times   \big\{   0    \big\}   \overset{\geq h}{\longleftrightarrow}  [    kn    , \big( k +2 \big) n ]   \times          \big\{    n  \big\}  \bigg]  >0 \text{ } \text{ , } \tag{\textbf{I}}
\end{align*}

\noindent or that,

\begin{align*}
  \textbf{P}_{h_2} \bigg[ \big(  kn    , \big( k + 1 \big) n ]   \times   \big\{ 0      \big\}      \overset{\geq h}{\longleftrightarrow}   [    kn    , \big( k + 1 \big) n ]  \times \big\{  n         \big\}  \bigg]  >0   \text{ } \text{ . } \tag{\textbf{II}}
\end{align*}

\noindent Both $(\textbf{I})$ and $(\textbf{II})$ are bound below by $\textbf{P}_{h_2} \big[ \mathscr{V}\mathscr{C} \big( n , 2n \big) \big]$. Furthermore,

\begin{align*}
       \textbf{P}_{h_2} \big[    \mathscr{H} \mathscr{C} \big( n , 2 N \big)        \big]    \geq \textbf{P}_{h_2} \big[    \mathscr{H} \mathscr{C} \big( n , 2 N \big)  \geq c \big( n , 2 N \big)  \big]  \geq \big(   1 - \textbf{P}_{h_2} \big[    \mathscr{V} \mathscr{C} \big( n , 2 N \big)     \big]     \big)^{2 \frac{N}{n}}             \text{ } \text{ , }
\end{align*}

\noindent because,

\begin{align*}
       1 - \textbf{P}_{h_2} \big[  \mathscr{V} \mathscr{C} \big( n , 2 N \big)   \big]  \geq      \big(   1 - \textbf{P}_{h_2} \big[    \mathscr{V} \mathscr{C} \big( n , 2 N \big)     \big]     \big)^{2 \frac{N}{n}}      \text{ } \text{ , }
\end{align*}

\noindent for the monotonic decreasing transformation,

\begin{align*}
       f \big( x \big) = \big( 1 - x \big)^{2 \frac{N}{n}}            \text{ } \text{ , }
\end{align*}

\noindent given $n,N$ satisfying,

\begin{align*}
   \frac{N}{n} < \frac{1}{2}  \text{ } \text{ . } 
\end{align*}

\noindent Next, observe that the \textit{horizontal crossing} between $N$ and $2N$ satisfies,

\begin{align*}
   \textbf{P}_{h_2} \big[    \mathscr{H} \mathscr{C} \big( n , 2 N \big)        \big]  \geq  \frac{N}{n} \big( 1 - \textbf{P}_{h_2}  \big[ \mathscr{V}\mathscr{C} \big( n , 2 N\big) ] \big)  \geq \frac{N}{n}   \big( 1 - \textbf{P}_{h_2}  \big[ \mathscr{V}\mathscr{C} \big( n , 2 N\big) ]  \big)^{2 \frac{N}{n}}  \\ \geq   \lfloor \frac{N}{n} \rfloor    \big( 1 - \textbf{P}_{h_2}  \big[ \mathscr{V}\mathscr{C} \big( n , 2 N\big) ]  \big)^{2 \frac{N}{n}} \\ \geq  \lfloor \frac{N}{n} \rfloor    \big( 1 - \textbf{P}_{h_2}  \big[ \mathscr{V}\mathscr{C} \big( n , 2 n \big) ]  \big)^{2 \frac{N}{n}}  \text{ } \text{ . }       
\end{align*}

\noindent Altogether,

\begin{align*}
                  \textbf{P}_{h_2} \big[    \mathscr{V}\mathscr{C} \big( n , 2n \big)   \big]     \leq \mathrm{exp} \bigg(  - \big( h_1 - h_2 \big)   \frac{N}{n} 
 \textbf{P}_{h_2} \big[ \mathscr{V} \mathscr{C} \big( n , 2n \big)  \big] \bigg)  \leq \mathrm{exp} \bigg(           - \big( h_1 - h_2 \big)   \frac{N}{n} \big(           1 - \textbf{P}_{h_2} \big[ \mathscr{V}\mathscr{C} \big( n , 2n \big) ]                    \big)^{2 \frac{N}{n}}                            \bigg) \text{ } \text{ , }
\end{align*}

\noindent from which we conclude the argument. \boxed{}

\bigskip

\noindent With the arguments below, we implement a similar inductive version for establishing that the first \textbf{Proposition} holds.

\bigskip

\noindent \textit{Proof of Proposition 1}. Define,

\begin{align*}
        \delta_k = \sqrt{\delta_{k+1}} \text{ } \text{ , } \\     n_k    =   n_{k+1}  \delta^2_k           \text{ } \text{ , }  \\   h_k     =   h_{k+1} + \delta_k \text{ } \text{ , } 
\end{align*}

\noindent recursively for each $k \geq 0$. From previous arguments, the exponential upper bound to the \textit{vertical crossing} would take the form,

\begin{align*}
        \textbf{P}_{h_{k+1}} \big[ \mathscr{V} \mathscr{C} \big( n_{k+1} , 2 n_{k+1} \big) \big]  \leq \mathrm{exp} \bigg(    - \big( h_k - h_{k+1} \big) \frac{n_{k+1}}{n_k} \big( 1 - \textbf{P}_{h_k} \big[     \mathscr{V} \mathscr{C} \big( n_{k+1} , 2 n_{k+1} \big)         \big]       \big)^{2 \frac{n_{k+1}}{n_k}}    \bigg)                \text{ } \text{ , } 
\end{align*}

\noindent which we can further manipulate to show,

\begin{align*}
         \textbf{P}_{h_{k+1}} \big[ \mathscr{V} \mathscr{C} \big( n_{k+1} , 2 n_{k+1} \big) \big]  \leq     \Delta_k                \text{ } \text{ , } 
\end{align*}

\noindent because,

\begin{align*}
     - \big( h_k - h_{k+1} \big) \frac{n_{k+1}}{n_k} \big( 1 - \textbf{P}_{h_k} \big[     \mathscr{V} \mathscr{C} \big( n_{k+1} , 2 n_{k+1} \big)         \big]       \big)^{2 \frac{n_{k+1}}{n_k}}    \leq   - \big( h_k - h_{k+1} \big) \frac{n_{k+1}}{n_k}   \big( 1 - \delta_k \big)^{2 \frac{n_{k+1}}{n_k}}         \leq  C \big( h_k , h_{k+1} \big)          \mathrm{log} \big(            \delta_k                 \big)  \\ \leq    \mathrm{log} \big( \Delta_k  \big)       \text{ } \text{ , } 
\end{align*}

\noindent for $\alpha$ sufficiently large, and parameters satisfying,

\begin{align*}
      \textbf{P}_{h_k} \big[     \mathscr{V} \mathscr{C} \big( n_{k+1} , 2 n_{k+1} \big)         \big]       \leq  \delta_k          \text{ } \text{ , } \\ - \big( h_k - h_{k+1} \big) \frac{n_{k+1}}{n_k} \leq   C \big( h_k , h_{k+1} \big)             \text{ } \text{ , } \\  \big( \mathrm{log} \big( n_{k+1} \big) - \mathrm{log} \big( n_k \big) \big) 2 \frac{n_{k+1}}{n_k} \big( 1 - \delta_k \big)      \leq   \delta_k       \text{ } \text{ , } \\   \frac{ \delta_k   }{\Delta_k} \leq 1    \text{ } \text{ . } 
\end{align*}

\noindent We conclude the argument by observing,

\begin{align*}
         \textbf{P}_{h_k - \epsilon} \big[  \mathscr{V}\mathscr{C} \big( N , 2 N \big)    \big]    \leq  \textbf{P}_{h_k - \epsilon} \big[  \mathscr{V}\mathscr{C} \big( n_k , 2 N \big)    \big]                   \overset{(\mathrm{*})}{\leq}    \big( \frac{n_0}{N} \big)^{\alpha - \frac{1}{2}}   \text{ } \text{ , } 
\end{align*}

\noindent where in $(\mathrm{*})$, we made use of the fact that,

\begin{align*}
     \frac{n_k}{2 n_{k+1}}  \big( \textbf{P}_{h_k - \epsilon} \big[  \mathscr{V}\mathscr{C} \big( n_k , 2 N \big)    \big] \big)^{\big( \alpha - \frac{1}{2} \big)^{-1}} \leq                     \overset{k}{\underset{i=0}{\prod}}   \delta^2_i       \leq    2 \frac{n_0}{ n_{k+1}} \leq   \frac{n_0}{N} \text{ } \text{ , } 
\end{align*}

\noindent for $\epsilon$ sufficiently small, implying,

\begin{align*}
  \textbf{P}_{h_k} \big[    0 \overset{\geq h}{\longleftrightarrow} \partial \Lambda_n    \big]   \leq                   2 \big(   \frac{n_k}{N}     \big)^{\frac{n_k}{n_{k+1}} - \frac{1}{2}}          \leq     2 \big(   \frac{n_k}{N}     \big)^{\alpha - \frac{1}{2}}                         \leq \mathrm{exp} \big( - c \big| x - y \big| \big)        \text{ } \text{ , } 
\end{align*}

\noindent for $c$ suitably large and $\alpha > 0$. Hence,

\begin{align*}
\underset{n \longrightarrow + \infty }{\mathrm{lim\text{ }  inf}}   \text{ }   \textbf{P} \big( \mathscr{V}   \mathscr{C} \big( n , 2 n \big)    \big)        \longrightarrow 0  \Rightarrow       \textbf{P}_{h_k} \big[    0 \overset{\geq h}{\longleftrightarrow} \partial \Lambda_n    \big]    \leq \mathrm{exp} \big( - c \big| x - y \big| \big)    \text{ } \text{ . } \boxed{} 
\end{align*}

\section{Crossing probabilities in the hard direction}

\noindent To control crossing probabilities in the hard direction, in comparison to arguments in the previous section for the easy direction, we concentrate on the following item.

\bigskip

\noindent \textbf{Proposition} \textit{2} (\textit{limit infimum of vertical crossings in the hard direction}). If $h \in \big( 0 , + \infty \big)$, there exists an infinite volume measure for which,

\begin{align*}
   \underset{n \longrightarrow + \infty}{\mathrm{lim \text{ } inf} }  \text{ }  \textbf{P}_h \big[ \mathscr{V} \mathscr{C} \big( n , 2n \big) ] > 0    \text{ } \text{ , } 
\end{align*}

\noindent for a \textit{vertical crossing}, then for any $h_0 > h$,

\begin{align*}
        \underset{n \longrightarrow + \infty}{\mathrm{lim \text{ } inf} }  \text{ }  \textbf{P}_{h_0} \big[ \mathscr{V} \mathscr{C} \big( n , 2n \big) ] > 0    \text{ } \text{ . } 
\end{align*}

\noindent To establish that the item above holds, introduce the item below.

\bigskip

\noindent \textbf{Lemma} \textit{3} (\textit{separated vertical crossings}). For $h \in \big( 0 , + \infty \big)$ and natural $n$, there exists an integer $I$, with $1 \leq I \leq \frac{n}{n^{\prime}}$, for $n^{\prime}$ sufficiently large, and strictly positive $c_0,c_1$, such that,

\begin{align*}
    I^2 \leq                   c_0 \text{ }  \frac{\textbf{P}_h \big[          \mathscr{V}\mathscr{C} \big( n , 2n \big)         \big]  }{\big( \textbf{P}_h \big[           \mathscr{V}\mathscr{C} \big( n , 2n \big)        \big] \big)^{\frac{c_1}{I}} }   \text{ } \text{ . } 
\end{align*}

\noindent This implies, for another strictly positive $c_3$,

\begin{align*}
  \textbf{P}_h \bigg[   \big|     \mathscr{V}\mathscr{C} \big(     \big[ 0 , 2n \big] \times \big[ 0 , \frac{n}{2} \big]  \big)             \big| =       2^I   \bigg]         \geq c_3 \textbf{P}_h \big[ \mathscr{V}\mathscr{C} \big( n , 2n \big) ]              \text{ } \text{ , } 
\end{align*}

\noindent where $\mathscr{V}\mathscr{C} \big(     \big[ 0 , 2n \big] \times \big[ 0 , \frac{n}{2} \big]     \big)$ denotes the \textit{vertical crossings} between $ \big[ 0 , 2n \big]$ and $\big[ 0 , \frac{n}{2} \big]$.

\bigskip

\noindent To prove the item above, we introduce the statement below.

\bigskip

\noindent \textbf{Corollary} \textit{1} (\textit{crossing probabilities in the hard direction are bound above by a doubly exponential function}). For some $\delta > 0$, there exists a strictly positive $c^{\prime\prime} \equiv c^{\prime\prime} \big( \delta \big)$, and $c^{\prime\prime\prime} \equiv c^{\prime\prime\prime} \big( \delta \big)$, such that for $h_2 > h_1$,

\begin{align*}
      \textbf{P}_{h_1} \big[        \mathscr{H}\mathscr{C} \big( \frac{n}{2} , n \big)          \big]    \leq \mathrm{exp} \bigg(        - c_3 \big( h_2 - h_1 \big)   \delta \text{ } \mathrm{exp} \big(                c_3                  \mathscr{I}   \big)     \bigg)          \text{ } \text{ , } 
\end{align*}

\noindent for,

\begin{align*}
      \mathscr{I} \equiv         f \big( \textbf{P}_{h_2} \big[   \mathscr{H}\mathscr{C} \big( n , 2n \big)         \big] \big)       \text{ } \text{ , } 
\end{align*}

\noindent where,

\[
   f \big( x \big) \equiv           \text{ }  
\left\{\!\begin{array}{ll@{}>{{}}l}    \frac{\mathrm{log} ( - {x}^{-1} ) }{\mathrm{log} ( \mathrm{log} ( - x^{-1} ) ) }       \text{ } \text{ , } \text{ } \text{for } x \in \big( -\infty , 0 \big) \text{ , }\\
       - \infty    \text{ } \text{otherwise}  \text{ . } \\
\end{array}\right.
\]

\noindent \textit{Proof of Corollary 1}. Fix all parameters as given in the statement above. In order to demonstrate that the doubly exponential upper bound holds, observe that there exists some upper bound for $I^2$, of the form,

\begin{align*}
     c_0 \text{ }  \frac{\textbf{P}_h \big[          \mathscr{V}\mathscr{C} \big( n , 2n \big)         \big]  }{\big( \textbf{P}_h \big[           \mathscr{V}\mathscr{C} \big( n , 2n \big)        \big] \big)^{\frac{c_1}{I}} }       \text{ } \text{ , } 
\end{align*}

\noindent for every $n \geq 1$, with $I$ defined by, 

\begin{align*}
    I \equiv \lfloor              c_{\textbf{f}}  \text{ }    \textbf{P}_{h_2} \big[  \mathscr{H} \mathscr{C} \big( n , 2n \big) \big]           \rfloor       \text{ } \text{ , } 
\end{align*}

\noindent for $c_{\textbf{f}}$ satisfying,

\begin{align*}
 \sqrt{\lfloor c_{\textbf{f}}  \rfloor }  \leq c_0                           \text{ } \text{ . } 
\end{align*}

\noindent Under the assumptions of \textbf{Lemma} \textit{3},

\begin{align*}
  \textbf{P}_h \bigg[   \big|     \mathscr{V}\mathscr{C} \big(     \big[ 0 , 2n \big] \times \big[ 0 , \frac{n}{2} \big]  \big)             \big| =       2^I   \bigg]         \geq  c \big( h \big) \textbf{P}_{h} \big(    \mathscr{V}\mathscr{C} \big(    n    ,   2n   \big)    \big) \geq c \big( h \big)  \frac{\delta}{2} \geq \frac{\delta}{2} \text{ } \text{ , }    
  \end{align*}

\noindent given $c \big( h \big)$ for which,

\begin{align*}
             c \big( h \big) \leq 1                 \text{ } \text{ . } 
\end{align*}

\noindent For $h_2 > h_1$,

\begin{align*}
    \frac{\textbf{P}_{h_2} \big(    \mathscr{H} \mathscr{C} \big( \frac{n}{2} , 2n       \big)     \big) }{\textbf{P}_{h_1} \big(       \mathscr{H} \mathscr{C} \big( \frac{n}{2} , 2n       \big)        \big) }    \leq     \mathrm{exp} \bigg(  - \big( h_2 - h_1 \big) \big( h_1 \big( 1 - h_1 \big)    \big)^{-1}     2^{I-2} \delta \bigg) \sim \text{ }           1  -  \bigg( \big( h_2 - h_1 \big) \big( h_1 \big( 1 - h_1 \big)    \big)^{-1}   2^{I-2}  \bigg) \delta     \text{ } \text{ . }
\end{align*}

\noindent Moreover,

\begin{align*}
      \mathrm{exp} \bigg(  - \big( h_2 - h_1 \big) \big( h_1 \big( 1 - h_1 \big)    \big)^{-1}     2^{I-2} \delta \bigg) \leq \mathrm{exp} \bigg(   -    \big( h_2 - h_1 \big)     \delta \mathrm{exp} \big(     c_3 \big( c_{\textbf{f}} \big) f \big(     \mathscr{H}\mathscr{C} \big(  n , 2n \big)    \big)    \big)           \bigg)       \text{ } \text{ , }
\end{align*}

\noindent for,

\begin{align*}
   \big( h_1 \big( 1 - h_1 \big)    \big)^{-1}  \leq 1         \text{ } \text{ , }
\end{align*}

\noindent $c_3 > 0$, and,

\begin{align*}
     I \leq \mathrm{log_2} \bigg(   \mathrm{exp}  \bigg(   f \big( \mathscr{H}\mathscr{C} \big( n , 2 n \big) \big)      \bigg)        \bigg)   \text{ } \text{ . }
\end{align*}

\noindent As a result,

\begin{align*}
      \textbf{P}_{h_2} \bigg[  \mathscr{H} \mathscr{C} \big( \frac{n}{2}, \frac{3}{5} n \big)      \cap  \mathscr{H} \mathscr{C} \big( \frac{3}{5} n,      \frac{7}{10} n  \big)  \cap  \mathscr{H} \mathscr{C} \big(  \frac{7}{10} n ,   \frac{4}{5}  n   \big)   \cap \mathscr{H} \mathscr{C} \big( \frac{4}{5} n  , \frac{9}{10} n  \big)  \cap      \mathscr{H} \mathscr{C} \big( \frac{9}{10} n , n  \big)     \bigg]    \overset{(\mathrm{FKG})}{\geq}  \underset{\mathscr{E} \in \mathscr{H}\mathscr{C}}{\mathrm{min}}  \text{ } \textbf{P}_{h_2} \big[         \mathscr{E}         \big]^5                     \text{ } \text{ , } 
\end{align*}

\noindent for the collection of events,

\begin{align*}
 \mathscr{H}\mathscr{C}  \equiv  \big\{  \mathscr{H} \mathscr{C} ( \frac{n}{2}, \frac{3}{5} n ) ,   \mathscr{H} \mathscr{C} ( \frac{3}{5} n ,      \frac{7}{10} n  )  ,  \mathscr{H} \mathscr{C} (  \frac{7}{10} n  ,   \frac{4}{5} n     )  ,  \mathscr{H} \mathscr{C} ( \frac{4}{5} n , \frac{9}{10} n )  , \mathscr{H} \mathscr{C} ( \frac{9}{10} n  , n  )   \big\}           \text{ } \text{ . } 
\end{align*}

\noindent The infimum obtained after applying $(\mathrm{FKG})$ can be upper bounded with,

\begin{align*}
      \textbf{P}_{h_2} \big[         \mathscr{H}\mathscr{C} \big( n , 2n \big)    \big]    \text{ } \text{ . }
\end{align*}

\noindent Concluding,

\begin{align*}
         \textbf{P}_{h_2} \big[         \mathscr{H}\mathscr{C} \big( \frac{n}{2} , n \big)    \big]  \leq   \mathrm{exp} \bigg(            -    \big( h_2 - h_1 \big)     \delta \mathrm{exp} \big(     c_3 \big( c_{\textbf{f}} \big) f \big(     \mathscr{H}\mathscr{C} \big(  n , 2n \big)    \big)    \big)           \bigg)     \text{ } \text{ . }     \boxed{} 
\end{align*}

\noindent With \textbf{Corollary} \textit{1}, below we provide arguments for \textbf{Proposition} \textit{2}.

\bigskip

\noindent \textit{Proof of Proposition 2}. Suppose $ \underset{n \geq 0 }{\mathrm{ \text{ } inf} }  \text{ }  \textbf{P}_h \big[ \mathscr{V} \mathscr{C} \big( n , 2n \big) ] > 0 $. With $\delta$, $c_0$ and $c_1$, introduce, for $h_0 > h$,

\begin{align*}
        n_k = 2^{-k} n_0    \text{ } \text{ , }  \\            h_k = h_0 - \big( h_0 - h \big) \sum_{i=1}^k 2^{-i}          \text{ } \text{ , } \\    \beta_k = \textbf{P}_{h_k} \big[         \mathscr{H}\mathscr{C} \big( n_k , 2 n_k \big)               \big]      \text{ } \text{ . }
\end{align*}

\noindent Next, consider,

\begin{align*}
 \textbf{P}_{h_2} \bigg[ \big\{ 0 \big\} \times [0 , n_k] \underset{[0 , \frac{n_k}{2} ] \times [0 , n_k ]}{\overset{\geq h}{\longleftrightarrow}} \big\{ \frac{n_k}{2} \big\} \times [0 , n_k]  \bigg]   \text{ } \text{ , } 
\end{align*}

\noindent which can be lower bounded, upon observing that,

\begin{align*}
            \textbf{P}_{h_2} \bigg[   \underset{i \in \textbf{Z}}{\bigcap}  \big\{   \big\{ 0 \big\} \times [ 0 , n_i ]  \underset{[0 , \frac{n_k}{2} ] \times [0 , n_k ]}{\overset{\geq h}{\longleftrightarrow}}     \big\{ \frac{n_i}{2} \big\} \times [ 0 , n_i ]                    \big\}     \bigg]   
 \overset{(\mathrm{FKG})}{\geq}     \underset{i \in \textbf{Z}}{\prod}   \textbf{P}_{h_2} \big[      \big\{ 0 \big\} \times [ 0 , n_i ]  \underset{[0 , \frac{n_k}{2} ] \times [0 , n_k ]}{\overset{\geq h}{\longleftrightarrow}}     \big\{ \frac{n_i}{2} \big\} \times [ 0 , n_i ]                       \big]  \geq C   C_1   \text{ , } 
\end{align*}

\noindent for some $C >0$,

\begin{align*}
 C_1 \equiv \textbf{P}_{h_2} \big[  \mathscr{H}\mathscr{C} \big(        \frac{n_k}{2}     ,   n_k    \big)    \big]   \text{ } \text{ , }  \tag{\textit{*}}
\end{align*}

\noindent and,

\begin{align*}
   C \equiv \underset{i \in \textbf{Z}}{\mathrm{inf}}  \text{ }    \textbf{P}_{h_2} \bigg[    \big\{ 0 \big\} \times [ 0 , n_i ]  \underset{[0 , \frac{n_k}{2} ] \times [0 , n_k ]}{\overset{\geq h}{\longleftrightarrow}}     \big\{ \frac{n_i}{2} \big\} \times [ 0 , n_i ]                       \bigg]    \text{ } \text{ . }
\end{align*}
\noindent Hence,

\begin{align*}
    \textbf{P}_{h_2} \big[           \mathscr{H}\mathscr{C} \big( n_k , 2 n_k \big)        \big]  \geq \frac{1}{C^2} \textbf{P}_{h_2} \big[              \mathscr{H}\mathscr{C} \big( n_k , 2 n_k \big)     \big]  \geq    \frac{C^4_1}{C^2}      \text{ } \text{ , } 
\end{align*}

\noindent as a result of the fact that,

\begin{align*}
 \underset{i \in \textbf{Z}}{\prod}   \textbf{P}_{h_2} \bigg[     \big\{ 0 \big\} \times [ 0 , n_i ]  \underset{[0 , \frac{n_k}{2} ] \times [0 , n_k ]}{\overset{\geq h}{\longleftrightarrow}}     \big\{ \frac{n_i}{2} \big\} \times [ 0 , n_i ]                        \bigg]  \geq \textbf{P}_{h_2} \big[              \mathscr{H}\mathscr{C} \big( n_k , 2 n_k \big)     \big]   \text{ } \text{ . } 
\end{align*}

\noindent Proceeding, to show that there exists $n_k$ for which,

\begin{align*}
     n_k > c^{\prime\prime} f \big( \textbf{P}_{h_2} \big[     \beta_k     \big] \big)   \text{ } \text{ , } 
\end{align*}

\noindent in light of the constant $C_1$ obtained in (\textit{*}), observe,

\begin{align*}
        \textbf{P}_{h_2} \big[            \mathscr{H}\mathscr{C} \big( n_k , 2 n_k \big)    \big]   \leq \frac{1}{\sqrt{n_k}} \equiv \frac{1}{\sqrt{2^{-k} n_0}} \leq \frac{1}{\sqrt{n_0}}                    \text{ } \text{ . } 
\end{align*}

\noindent However,

\begin{align*}
           \textbf{P}_{h_2} \big[            \mathscr{H}\mathscr{C} \big( n_k , 2 n_k \big)    \big] \geq     \delta       \text{ } \text{ , } 
\end{align*}

\noindent contradicts the fact that the following upper bound holds,

\begin{align*}
    \textbf{P}_{h_2} \big[             \mathscr{H}\mathscr{C} \big( n_k , 2 n_k \big)      \big] \leq     \textbf{P}_{h_k} \big[                \mathscr{H}\mathscr{C} \big( n_k , 2 n_k \big)         \big]     \equiv \beta_k   \leq   \mathrm{exp} \big(  - \frac{n_k}{c^{\prime\prime}}    \big)  \leq \frac{1}{n^l_k}      \text{ } \text{ , }
\end{align*}

\noindent for $l$ sufficiently large. From previous arguments, the inequality in terms of $\beta_k$ and $\beta_{k+1}$,

\begin{align*}
       \beta_{k+1}    \leq   \mathrm{exp} \bigg( -  c^{\prime\prime\prime}         \frac{n_k}{n_{k+1}}    \big( h_2 - h_1 \big)    \delta \mathrm{exp} \big( c^{\prime\prime\prime}     f \big( \beta_k \big)       \big)     \bigg)       \text{ } \text{ , } 
\end{align*}

\noindent for strictly positive $c^{\prime\prime\prime} \equiv c^{\prime\prime\prime} \big( \delta \big)$, implies the existence of a constant $c^{\prime\prime\prime\prime}$ for which,

\begin{align*}
       - c^{\prime\prime\prime} \frac{n_k}{n_{k+1}} \big( h_2 - h_1 \big)          \mathrm{exp} \big( c^{\prime\prime\prime}     f \big( \beta_k \big)       \big)  \leq   - c^{\prime\prime\prime}  \big( h_2 - h_1 \big)          \mathrm{exp} \big( c^{\prime\prime\prime\prime}     f \big( \beta_k \big)   \big)        \text{ } \text{ . } 
\end{align*}

\noindent Hence,

\begin{align*}
          \beta_{k+1}     \leq   \mathrm{exp} \bigg( -    c^{\prime\prime\prime}  \big( h_2 - h_1 \big)          \mathrm{exp} \big( c^{\prime\prime\prime\prime}     f \big( \beta_k \big)   \big)  \bigg)  \sim 1 -   c^{\prime\prime\prime}  \big( h_2 - h_1 \big)          \mathrm{exp} \big( c^{\prime\prime\prime\prime}     f \big( \beta_k \big)   \big) \beta_k \leq 1 -   \beta_k \Delta^{-1} \leq       c_{\Delta} \Delta^{k-1}           \text{ } \text{ , } 
\end{align*}

\noindent for $\Delta^{-1} \geq \mathrm{exp} \big(   c^{\prime\prime\prime}     f \big( \beta_k \big)    \mathrm{log} \big(  c^{\prime\prime\prime\prime}  \big( h_2 - h_1 \big)   \big)      \big)$, and $k > 0$ for which,

\begin{align*}
       \mathrm{log} \big(  \frac{ \Delta  - \beta_k }{c_{\Delta}} \big)  \approx    k     \text{ } \text{ . } 
\end{align*}

\noindent Suppose $\beta_0 \leq c_{\Delta}$. For any $k$ such that $k \geq \mathrm{log}_2 \big( \sqrt{n_0} \big)$, there exists $m \in [ \sqrt{n_0} , n_0 ]$ for which,

\begin{align*}
       \textbf{P}_{h_2} \big[   \mathscr{H} \mathscr{C} \big( m , 2m \big)      \big] \leq \beta_0 \text{ } \mathrm{exp} \big( -   c^{\prime\prime\prime}    f \big( \beta_k \big)    \mathrm{log} \big(  c^{\prime\prime\prime\prime}  \big( h_2 - h_1 \big)   \big)                \big)  \leq c_{\Delta} \text{ }  \mathrm{exp} \big(         -  c^{\prime\prime\prime}   \lfloor   f \big( \beta_k \big)    \mathrm{log} \big(  c^{\prime\prime\prime\prime}  \big( h_2 - h_1 \big)    \big)  \rfloor        \big)   \leq c_{\Delta} \text{ }                   m^N     \text{ } \text{ . } 
\end{align*}

\noindent for $N$ sufficiently large. As with previous computations, we obtain a contradiction again, as,

\begin{align*}
  \textbf{P}_{h_2} \big[ \mathscr{H}\mathscr{C} \big( m , 2m \big)  \big]  \leq   \textbf{P}_{h_k} \big[ \mathscr{H}\mathscr{C} \big( m , 2m \big)  \big] \equiv \beta_k \leq \frac{1}{n^{l^{\prime}}_k} \text{ } \text{ , } 
\end{align*}

\noindent for $l^{\prime}$ sufficiently large. Hence,

\begin{align*}
      \underset{n \longrightarrow + \infty}{\mathrm{lim \text{ } inf} } \text{ }  \textbf{P}_{h_1} \big[ \mathscr{H} \mathscr{C} \big( n , 2n \big) \big] \geq c_{\Delta} > 0     \text{ } \text{ . } \boxed{}
\end{align*}

\noindent We now provide arguments for $\textbf{Lemma}$ \textit{3}.

\bigskip

\noindent \textit{Proof of Lemma 3}. Fix all constants as given in the statement of \textbf{Lemma} \textit{3}. In order to combine crossings in the hard direction to obtain a crossing of $[0,2n] \times [0,n]$, define,

\begin{align*}
       \alpha \equiv   \underset{k \in [ \frac{n}{8}, \frac{n}{2}]}{\mathrm{sup}} \textbf{P}_{h_2} \big[   \mathscr{H} \mathscr{C} \big( \lceil  \big( 2 + v \big) k         \rceil , 2 k  \big)      \big] \text{ } \text{ , } 
\end{align*}

\noindent for $v \equiv 10^{-2} I$. To combine crossings in the hard direction to obtain crossings across $[0,2n] \times [0,n]$, write,

\begin{align*}
        \alpha \overset{(\mathrm{FKG})}{\leq} \underset{v \in \textbf{R} :  \lceil ( 2 + v ) k \rceil < 2 k    }{\prod} \textbf{P}_{h_2} \big[    \mathscr{H} \mathscr{C} \big( \lceil  \big( 2 + v \big) k         \rceil , 2 k  \big)    \big]   {\leq}   \textbf{P}_{h_2} \big[  \mathscr{H}\mathscr{C} \big( n , 2n \big)          \big]^{\frac{v}{N^{\prime}}}     \leq  \textbf{P}_{h_2} \big[   \mathscr{H}\mathscr{C} \big( n , 2n \big)         \big]^{\frac{1}{N^{\prime} I }} \\  \leq  \textbf{P}_{h_2} \big[   \mathscr{H}\mathscr{C} \big( n , 2n \big)         \big]^{\frac{c^{\prime}}{ I }}  \text{ } \text{ , } 
\end{align*}

\noindent for $N^{\prime} > 32$, and $c^{\prime} \geq \frac{1}{N^{\prime}}$. Altogether, to establish that the claim holds, one must argue that,

\begin{align*}
        \textbf{P}_{h_2} \bigg[           \big\{  - \big( 1 + u \big) k  \big\} \times [ 0 , 2k ]     {\overset{\geq h}{\longleftrightarrow}}    \big\{  \big( 1 + u \big) k \big\} \times [ 0 , 2k ]     \bigg]              \text{ } \text{ , } 
\end{align*}

\noindent occurs with positive probability, which we denote as the event $\mathscr{C} \big( R \big( k \big) \big)$, as well as,

\begin{align*}
  \textbf{P}_{h_2} \bigg[ \forall v_i , v_j \in  \Lambda  :  v_i \cap \big( [ - 3 uk , 3 uk] \times \{ 0 \} \big) = \emptyset ,  v_j \cap \big( [ - 3 uk , 3 uk] \times \{ 2k \} \big) = \emptyset      \bigg]  \text{ } \text{ , } 
\end{align*}

\noindent occurring with vanishing probability, for vertices satisfying the condition,

\begin{align*}
   v_i , v_j \in \mathscr{C} \big( R \big( k \big) \big) \Longleftrightarrow  v_i \cap  \mathscr{C} \big( R \big( k \big) \big)  \neq \emptyset , v_j \cap  \mathscr{C} \big( R \big( k \big) \big)  \neq \emptyset    \text{ } \text{ , } 
\end{align*}

\noindent and the finite volume $\Lambda$, satisfying the condition,

\begin{align*}
    \textbf{P}_{h_2} \bigg[    \forall v \in \Lambda , \exists \text{ countably many } v_k \in \mathscr{C} \big( R \big( k \big) \big) :   v \cap v_k \neq \emptyset     \bigg]        \text{ } \text{ , } 
\end{align*}

\noindent which we denote as the event $\mathscr{E} \big( \mathscr{C} \big( R \big( k \big) \big) , k \big) \equiv \mathscr{E} \big(  k \big)$. To upper bound $ \textbf{P}_{h_2} \big[ \mathscr{E} \big( k \big) \big]$, observe,

\begin{align*}
        \beta - \alpha = \textbf{P}_{h_2} \big[            \mathscr{V}\mathscr{C}    \big]            \text{ } \text{ , } 
\end{align*}

\noindent for paths $\gamma \in \Gamma$ for which a \textit{vertical crossing} occurs,

\begin{align*}
   \mathscr{V}\mathscr{C} \equiv    \underset{u \in \textbf{Z}}{\bigcup}     \big\{  \mathrm{paths \text{ } } \gamma   \text{ }         \big|    \gamma \cap   \big(   [ - \big( 1 + u \big) k , - 3 u k ]    \times \{  0      \}      \big)          \neq \emptyset                \big\}      \text{ } \text{  , }
\end{align*}

\noindent from which one obtains,

\begin{align*}
        \textbf{P}_{h_2} \big[ \mathscr{V}\mathscr{C} \cap  \mathscr{V}\mathscr{C}^{\prime}    \big]   \overset{(\mathrm{FKG})}{\geq}  \textbf{P}_{h_2} \big[   \mathscr{V}\mathscr{C}  \big]     \textbf{P}_{h_2} \big[   \mathscr{V}\mathscr{C}^{\prime}  \big]   \geq \big( \beta - \alpha \big) \gamma            \text{ } \text{ , }
\end{align*}

\noindent where the \textit{vertical crossing} $\mathscr{V}\mathscr{C}^{\prime}$ is,

\begin{align*}
   \mathscr{V}\mathscr{C}^{\prime} \equiv    \underset{u \in \textbf{Z}}{\bigcup}     \big\{  \mathrm{paths \text{ } } \gamma   \text{ }         \big|    \gamma \cap   \big(   [ - \big( 1 + 4 u \big) k , - \big( 1 - 2 u \big) k ]    \times \{  0      \}      \big)          \neq \emptyset                \big\}      \text{ } \text{  . }
\end{align*}

\noindent On the other hand, if we consider a similar \textit{vertical crossing} $R^{\prime} \big( k \big)$ for which $R \big( k \big)$ occurs but does not intersect $ [ - \big( 1 + 4 u \big) k , \big( 1 - 2 u \big) k ] \times [ 0, 2k]$, the upper bound instead takes the form,

\begin{align*}
 \textbf{P}_{h_2} \big[   R^{\prime} \big( k \big)               \big] \equiv \textbf{P}_{h_2} \bigg[       R \big( k \big) \cap \bigg(       R \big( k \big) \overset{\geq h}{\not\longleftrightarrow}                   \big(  [ \big( 1 + 4 u \big) k , \big( 1 - 2 u \big) k ] \times [ 0, 2k]        \big)            \bigg)   \bigg]  \leq       \alpha + \alpha^{\prime}       \text{ } \text{ , }       \tag{Bound 1}
\end{align*}

\noindent for $0 < \alpha^{\prime} < 1$. The second upper bound estimate asserts,

\begin{align*}
    \textbf{P}_{h_2} \big[ R^{\prime\prime} \big( x \big)  \big]   \equiv \textbf{P}_{h_2} \bigg[ R \big( x \big) \cap             \big( R \big( x \big) \overset{\geq h}{\not\longleftrightarrow}    \big( \big\{ \big( 1 - 2u \big) k \big\} \times \big[ 0 , 2k \big] \big)         \big)  \bigg]           \leq \alpha + \gamma^{\prime\prime}    \text{ } \text{ , }       \tag{Bound 2}
\end{align*}

\noindent where $\gamma$ is a path such that $R \big( x \big)$ occurs, and sufficient $\gamma^{\prime\prime}$.

\bigskip

\noindent The third upper bound estimate asserts,

\begin{align*}
         \textbf{P}_{h_2} \big[     R^{\prime\prime\prime} \big( x \big)      \big] \equiv  \textbf{P}_{h_2} \big[    \gamma_1 \big(   \textbf{R}   ,  \big[ 0 , \big( 2 - 11 u \big) k \big]   \big)   \cap             \mathscr{HC} \big( \gamma_1 , \gamma_2 \big)                  \cap  \gamma_2 \big(   \big( 1 - 2u \big) k       ,   \big[ 0 , \big( 2 - 11 u \big) k \big]   \big)                    \big]             \text{ } \text{ , } 
\end{align*}

\noindent where $\gamma_1$ and $\gamma_2$ denote paths, the first of which is, 

\begin{align*}
     \gamma_1 \equiv \gamma_1 \big( \textbf{R}  ,   \big[ 0 , \big( 2 - 11 u \big) k \big]  \big) \equiv  \textbf{R} \times   \big[ 0 , \big( 2 - 11 u \big) k \big]         \text{ } \text{ , } 
\end{align*}

\noindent the second of which is,

\begin{align*}
      \gamma_2 \equiv  \gamma_2 \big(  \big( 1 - 2u \big) k  , \big[ 0 , \big( 2 - 11 u \big) k \big]  \big) \equiv  \big( 1 - 2 u \big) k      \times  \big[ 0 , \big( 2 - 11 u \big) k \big]    \text{ } \text{ , } 
\end{align*}

\noindent and,

\begin{align*}
   \mathscr{H}\mathscr{C} \big( \gamma_1 , \gamma_2 \big) \equiv \big\{   [ - 3 u k , 3 uk] \times \big\{ 0 \big\}     \overset{\geq h}{\longleftrightarrow}  \big\{ \big( 1 - 2 u \big) k \big\} \times \big[ 0 , \big( 2- 11 u \big) k \big]       \big\}   \text{ } \text{ . } 
\end{align*}

\noindent Depending upon whether paths for which $\mathscr{H}\mathscr{C} \big( \gamma_1 , \gamma_2 \big)$ occurs intersect $ \big\{    - \big( 1 - 8 u \big) k \big\} \times \big[ 0 , \big( 2 - 11 u \big) k \big]$, the probability $\textbf{P}_{h_2} \big[ R^{\prime\prime} \big( x \big) \big]$ admits the upper bound,

\begin{align*}
       \textbf{P}_{h_2} \big[ R^{\prime\prime\prime} \big( x \big) \big] \leq         \alpha +   \alpha^{\prime\prime}      \text{ } \text{ , }           \tag{Bound 3, I}
\end{align*}

\noindent while if the intersection with $ \big\{    - \big( 1 - 8 u \big) k \big\} \times \big[ 0 , \big( 2 - 11 u \big) k \big]$ does not occur, the probability $\textbf{P}_{h_2} \big[ R^{\prime\prime} \big( x \big) \big]$ admits the upper bound,

\begin{align*}
 \textbf{P}_{h_2} \big[ R^{\prime\prime\prime}_{\emptyset} \big( x \big)     \big]  \equiv  \textbf{P}_{h_2} \big[ R^{\prime\prime} \big( x \big) \big]  \leq                     \beta - \alpha +       \alpha^{\prime\prime}    \text{ } \text{ . } \tag{Bound 3, II}
\end{align*}

\noindent Comparing the upper bounds for the probability of $R^{\prime\prime} \big( x \big)$ occurring implies the following upper bound for each possible path,

\begin{align*}
    \beta - \alpha + \alpha^{\prime\prime} \leq \alpha + \alpha^{\prime\prime} \leq    \alpha + \gamma^{\prime\prime\prime}   \text{ } \text{ , } 
\end{align*}

\noindent for sufficient $\gamma^{\prime\prime\prime}$. Hence,

\begin{align*}  
     \textbf{P}_{h_2} \big[ R^{\prime\prime\prime}_{\emptyset} \big( x \big)     \big]  \leq \textbf{P}_{h_2} \big[ R^{\prime\prime\prime} \big( x \big)     \big]  \leq \alpha + \gamma^{\prime\prime\prime}  \text{ } \text{ . } \tag{Bound 3}
\end{align*}

\noindent The fourth upper bound is concerned with crossings of rectangles $\mathscr{R}_j$, where,

\begin{align*}
\mathscr{R} \equiv \big[ 0 , 2n  \big] \times \big[  - k , k    \big]  \supsetneq   \mathscr{R}_j \equiv                \underset{u \in \textbf{Z}}{\bigcup}  \big\{      \big[     j u k        ,      \big( 2 + \big( j + 2 \big) u \big) k       \big]              \times \big[  - k , k    \big]               \big\}     \text{ } \text{ , } 
\end{align*}

\noindent for $0 \leq j \leq J$, where,

\begin{align*}
          J =     \lfloor       \frac{1}{u} \big( \frac{n}{k} - 2 \big)          \rfloor  - 2      \text{ } \text{ . } 
\end{align*}

\noindent The upper bound is equivalent of $\big\{ R^{\prime} \big( k \big) \cap R^{\prime\prime}\big( k \big) \cap   R^{\prime\prime\prime} \big( x \big)   \big\}$ to,

\begin{align*}
     \textbf{P}_{h_2} \big[        R^{\prime} \big( k \big)  \cap            R^{\prime\prime} \big( k \big)           \cap          R^{\prime\prime\prime} \big( k \big)   \big] \overset{(\mathrm{Bound \text{ } 1}) , (\mathrm{Bound \text{ } 2}) , (\mathrm{Bound\text{ }  3})}{\leq} \frac{3}{u} \text{ }  \mathrm{max} \big\{ \alpha , \alpha^{\prime} , \gamma^{\prime\prime\prime} \big\}    \text{ } \text{ . } 
\end{align*}

\noindent Next, consider paths $\gamma^1 \cap \gamma^2 \neq \emptyset$ for which,

\begin{align*}
      \textbf{P}_{h_2} \big[  \mathrm{paths \text{ }} \gamma^1 :    \gamma^1 \cap  \mathcal{I}_{-k}    \neq \emptyset  \big]  > 0        \text{ } \text{ , } \\ 
       \textbf{P}_{h_2} \big[   \mathrm{paths \text{ }} \gamma^2 :    \gamma^2 \cap \mathcal{I}_k \neq \emptyset     \big]    > 0   \text{ } \text{ , }
\end{align*}

\noindent for,

\begin{align*}
      \mathcal{I}_{-k}   \equiv        \big[ \big( 1 + \big( j - 2 \big) u \big) k , \big( 1 + \big( j + 4 \big) u \big) k ]   \times \big\{ - k \big\}       \text{ } \text{ , } 
\end{align*}

\noindent and,

\begin{align*}
      \mathcal{I}_k    \equiv     \big[ \big( 1 + \big( j - 2 \big) u \big) k , \big( 1 + \big( j + 4 \big) u \big) k ]   \times  \big\{     k    \big\}  \text{ } \text{ . } 
\end{align*}

\noindent Between $\gamma^1$ and $\gamma^2$, the fact that $\big\{ R^{\prime} \big( x \big) \cap R^{\prime\prime} \big( x \big) \cap   R^{\prime\prime\prime} \big( x \big)   \big\}$ does not occur implies,

\begin{align*}
            \textbf{P}_{h_2} \big[          \mathscr{V} \mathscr{C} \big(      \mathscr{D} , \geq h \big)             \big]  = 0            \text{ } \text{ , } 
\end{align*}

\noindent where,

\begin{align*}
    \mathscr{D}   \equiv  \gamma^1 \cap \mathcal{I}_{-k} \cap \mathcal{I}_k \cap \gamma^2    \text{ } \text{ , } 
\end{align*}

\noindent and, $\mathscr{V}\mathscr{C} \big( \mathscr{D} , \geq h \big)$ denotes the \textit{vertical crossing} across $\mathscr{D}$ of height $\geq h$. The same argument holds for crossings across $\big[ 0 , 2n \big] \times [- k , k \big]$. Fix,

\begin{align*}
k_i \equiv      \lfloor   \big(  1 - n^{\prime} v i       \big)   \frac{n}{2}    \rfloor               \text{ } \text{  . }
\end{align*}

\noindent for some $n^{\prime} > n$. To conclude the proof of the lemma, observe that the intersection of \textit{horizontal crossings} can be upper bound with a \textit{vertical crossing}, as,

\begin{align*}
    \textbf{P}_{h_2} \bigg[  \text{ }        \overset{I-1}{\underset{i=0}{\bigcup}}    \mathscr{H}\mathscr{C} \big( k_i \big)      \bigg]   \leq   \frac{3}{u} \text{ }  \mathrm{max} \big\{ \alpha , \alpha^{\prime} , \gamma^{\prime\prime\prime} \big\}   \leq     300   I^2    \mathrm{max} \big\{ \alpha , \alpha^{\prime} , \gamma^{\prime\prime\prime} \big\}   \leq             \mathscr{C}   \mathrm{max} \big\{ \alpha , \alpha^{\prime} , \gamma^{\prime\prime\prime} \big\}  \leq \mathscr{C} \text{ } \textbf{P}_{h_2} \big[ \mathscr{C} \mathscr{V} \big( n , 2n \big) \big] \text{ } \text{ , } 
\end{align*}

\noindent for $I \geq \sqrt{\frac{100}{u}}$, and $\mathscr{C} < 1$. Hence,

\begin{align*}
 \textbf{P}_{h_2} \bigg[   \big| \mathscr{V} \mathscr{C} \big(      \big[ 0 , 2n \big]        \times \big[ - k , k  \big]        \big) \big| \equiv 2^I       \bigg]    \geq \textbf{P}_{h_2} \big[         \mathscr{C} \mathscr{V} \big( n , 2n \big)    \big]      \text{ } \text{ , } 
\end{align*}

\noindent in which the \textit{vertical crossings} across $\big[ 0 , 2n \big] \times \big[ - k_I , k_I \big]$ occur with at least probability $\textbf{P}_{h_2}\big[               \mathscr{V} \mathscr{C}  \big( n , 2n \big)     \big]$, as,

\begin{align*}
     \textbf{P}_{h_2} \bigg[        \big| \mathscr{V} \mathscr{C} \big(        \big[ 0 , 2n \big] \times \big[ - k_I , k_I \big]   \big)        \bigg]          \geq \textbf{P}_{h_2} \big[ \mathscr{V}\mathscr{C}  \big( n , 2n \big) \big]   \text{ } \text{ . } \boxed{}
\end{align*}

\bigskip

\noindent To prove the main theorem which will establish the sharpness of the phase transition, we implement the following arguments. We make use of the terms $\textbf{P} \big( A \big|   \phi_x \geq h      \big)$ and $ \textbf{P} \big( A \big|    \phi_x < h     \big)$ in the definition of the \textit{differential inequality} for the GFF. For the proof, we rely upon the following two statements.

\bigskip

\noindent \textbf{Proposition} \textit{1} (\textit{almost-sure dominance}). For a finite graph $G$, vertex $x$, and boundary condition $\xi$, there exists a product measure $\Phi$ on $\Omega \times \Omega$, which satisfies,

\begin{align*}
     \Phi \text{ almost surely } \pi \leq \omega \text{, for edges } f \not\in C_x \big( \omega \big)          \text{ } \text{ , } 
\end{align*}

\noindent for the height $h$ of the GFF, the \textit{open cluster} $C_x \big( \omega \big)$ about $x$, $\big( \pi , \omega \big) \sim \Phi$, where,

\begin{align*}
   \pi \sim \textbf{P}_h \big[ \cdot | \varphi_x < h \big] \text{ } \text{ , } \\    \omega \sim \textbf{P} \big[ \cdot | \varphi_x \geq h \big]    \text{ } \text{ . } 
\end{align*}

\bigskip

\noindent The next item below establishes how \textit{horizontal crossings} across the free field are related for an increasing sequence of \textit{height} parameters.

\bigskip

\noindent \textbf{Corollary} \textit{2} (\textit{coupling crossings across the Gaussian free field for different height parameters}). For any $ - \infty < h_0 < h_1 < 0$, there exists a strictly positive constant $c \equiv c \big( h_0 \big)$ so that,

\begin{align*}
    \textbf{P}_{h_0} \big[ \mathscr{H}\mathscr{C} \big(        n , 2n  \big)  \big] \big( 1 -    \textbf{P}_{h_1} \big[ \mathscr{H}\mathscr{C} \big(      n , 2n    \big)  \big]  \big) \leq \big(    \textbf{P}_{h_1} \big[ 0    \overset{\geq h}{\longleftrightarrow} \partial \Lambda_n  \big]  \big)^{c ( h_1 - h_0 )}     \text{ } \text{ , } 
\end{align*}

\noindent for any $n \geq 1$. 

\bigskip

\noindent \textit{Proof of Proposition 1}. Fix $G$, $e \in E_G$, $\chi$, $\Omega$, and $h$. Over the state space $\Omega \times \Omega$, for a function $f$ over $E_G$, denote the two configurations $\omega^f$ and $\omega_f$ as the configurations which are equal to $1$ and $0$, respectively. From $f$, denote the indicator $\textbf{1} \big( f , \omega \big)$ for the event that there exists endpoints of $f$ which are not connected in $\omega^f \backslash \{ f \}$.

\bigskip

\noindent Over $\Omega \times \Omega$, denote the continuous time Markov chain,

\begin{align*}
    S \equiv \big\{     \big( \pi , \omega \big) \in \Omega \times \Omega : \varphi_x  < h , \varphi_x \geq h, \pi \leq \omega , \pi \big( f \big) \neq \omega \big( f \big) \text{ } \text{ , } \text{ } \forall f \not\in C_e \big( \omega \big)    \big\}     \text{ } \text{ . } 
\end{align*}

\noindent From the generators $J$ for $S$, in a previous reference (see {\color{blue}[8]}), it has already been proven that the Markov chain has a unique invariant measure from the coupling $\Phi$. \boxed{}

\bigskip

\noindent \textit{Proof of Corollary 2}. Fix $h_0$ and $h_1$ as state in the Corollary. Under the assumption that such a unique infinite-volume measure exists, for $G^{\prime} \subset G$ such that $G \cap \big( \big[ 0 , 2n \big] \times \big[ 0 , n \big] \big) \neq \emptyset$, and $\forall v \in V_G$, denote,

\begin{align*}
     \textbf{P}_{G^{\prime},h } \big[ \cdot \big] \equiv \textbf{P}_h \big[ \cdot \big]   \text{ } \text{ , } 
\end{align*}

\noindent from which we write,

\begin{align*}
  \textbf{P}_h \big[ \mathscr{H} \mathscr{C} \big( n , 2n \big) \big| \varphi_x \geq h \big]  -  \textbf{P}_h \big[ \mathscr{H} \mathscr{C} \big( n , 2n \big) \big| \varphi_x < h  \big]          = \Phi \big[   w \in \mathscr{H}\mathscr{C} \big( n , 2n \big) :  \pi \not\in \mathscr{H}\mathscr{C} \big( n , 2n \big)      \big]                   \text{ } \text{ . }
\end{align*}

\noindent Under $\Phi \big[ \cdot \big]$, the event $\big\{     w \in \mathscr{H}\mathscr{C} \big( n , 2n \big) :  \pi \not\in \mathscr{H}\mathscr{C} \big( n , 2n \big)         \big\}$ occurring implies,

\begin{align*}
    \textbf{P}_{h} \big[  \mathscr{H}\mathscr{C} \big( n , 2n \big)         \big| \varphi_x \geq h  \big]  - \textbf{P}_h \big[ \mathscr{H}\mathscr{C} \big( n , 2n \big)     \big|   \varphi_x < h \big]   \leq \Phi \big[   u \underset{\omega , G}{\overset{\geq h}{\longleftrightarrow}} \Lambda_n + u   \big] + \Phi \big[  v    \underset{\omega , G}{\overset{\geq h}{\longleftrightarrow}}  \Lambda_n + u    \big] \\  \leq c^{\prime} \textbf{P}_h \big[ u \overset{\geq h}{\longleftrightarrow} \partial \Lambda_n + u  \big]  \text{ } \text{ , }
\end{align*}

\noindent for $c^{\prime}$ sufficiently large, where under this choice of $x$, the radius of the open cluster around $C_x \big( \omega \big)$ has radius at least $n$. Finally, to prove that the desired inequality holds, we make use of the \textit{differential inequality} for the GFF, in which,

\begin{align*}
\frac{\mathrm{d}}{\mathrm{d}h} \text{ } \mathrm{log} \bigg[       \frac{\textbf{P}_h \big[   \mathscr{H}\mathscr{C} \big( n , 2n \big)     \big]}{1 - \textbf{P}_h \big[   \mathscr{H}\mathscr{C} \big( n , 2n \big)             \big]}   \bigg] \equiv -  \frac{\mathrm{d}}{\mathrm{d}h} \text{ } \mathrm{log} \bigg[       \frac{1 - \textbf{P}_h \big[   \mathscr{H}\mathscr{C} \big( n , 2n \big)             \big]}{\textbf{P}_h \big[   \mathscr{H}\mathscr{C} \big( n , 2n \big)     \big] }   \bigg] \geq    - c  \text{ }  \mathrm{log} \big[        \underset{u \in V_G}{\mathrm{max}}  \textbf{P}_h \big[  u       \overset{\geq h}{\longleftrightarrow}  \partial \Lambda_n + u       \big]     \big]     \text{ }\text{ , } 
\end{align*}

\noindent for strictly positive $c$. Observing that the RHS of the inequality above is decreasing, the LHS can also be arranged as,

\begin{align*}
 \frac{\mathrm{d}}{\mathrm{d}h} \text{ } \bigg(   \mathrm{log} \big[    \textbf{P}_h \big[   \mathscr{H}\mathscr{C} \big( n , 2n \big)     \big]    \big] -      \mathrm{log} \big[     1 - \textbf{P}_h \big[   \mathscr{H}\mathscr{C} \big( n , 2n \big)             \big]   \big]       \bigg)  \text{ } \text{ . } 
\end{align*}

\noindent Integrating the inequality above between $h_0$ and $h_1$ yields, for one side of the differential inequality,

\begin{align*}
     \text{ } \int_{h=h_0}^{h_1} \mathrm{log} \big[        \underset{u \in V_G}{\mathrm{max}}  \textbf{P}_h \big[  u       \overset{\geq h}{\longleftrightarrow}  \partial \Lambda_n + u       \big]     \big]^{- c }   \mathrm{d} h   \equiv    \big(      \underset{u \in V_G}{\mathrm{max}}  \textbf{P}_h \big[  u       \overset{\geq h}{\longleftrightarrow}  \partial \Lambda_n + u       \big]        \big)^{c ( h_1-h_0)}   \text{ } \text{ , }
\end{align*}

\noindent which in turn yields, in combination with rearrangements of the other side of the differential inequality,

\begin{align*}
   \textbf{P}_{h_0} \big[    \mathscr{H}\mathscr{C} \big( n  , 2n \big)    \big]    \big(  1 -    \textbf{P}_{h_1}   \big[    \mathscr{H}\mathscr{C} \big( n  , 2n \big)    \big]      \big)    \leq     \big(      \underset{u \in V_G}{\mathrm{max}}  \textbf{P}_{h_1} \big[  u       \overset{\geq h}{\longleftrightarrow}  \partial \Lambda_n + u       \big]        \big)^{c ( h_1-h_0) }    \text{ } \text{ . }
\end{align*}

\noindent Taking the finite volume limit as $G^{\prime} \longrightarrow G$ yields the desired result. \boxed{}

\bigskip

\noindent Next, introduce the following item for exponential decay between two vertices on the dual graph.

\bigskip

\noindent \textbf{Proposition} \textit{3} (\textit{exponential decay in the dual graph}). For $h \in \big( - \infty , 0 \big)$, and the finite-volume measure $\textbf{P}_h \big[ \cdot \big]$, there exists $c \equiv c \big( h \big)$ such that,

\begin{align*}
\textbf{P}_h \big[   u \underset{*}{\overset{\geq h}{\longleftrightarrow}} v   \big] \leq \mathrm{exp} \big( - c \big| u - v \big| \big) \text{ } \text{ , } 
\end{align*}

\noindent over the dual graph $G^{*}$ to $G$. As a result, the probability of obtaining an infinite connected component vanishes, in which $\textbf{P}_h \big[ 0 \overset{\geq h}{\longleftrightarrow} + \infty \big] = 0$, and $h \geq h_c$.

\bigskip

\noindent \textit{Proof of Proposition 3}. Fix $u ,v \in G^{*}$. To demonstrate that the event $\big\{  u \underset{*}{\overset{\geq h}{\longleftrightarrow}} v  \big\}$ occurs with exponentially small probability proportional to $-c \big| u - v \big|$, consider the dual crossing event across the annulus $A \big( v \big)$. For the dual crossing event to occur, there must exist an open path surrounding $0 \in G^{*}$ intersecting $v$. Under the assumption that the dual graph is locally finite, the possible number of vertices $u$ for which such a dual path exists is bound from above by $C \big| v \big|$, for finite $C \equiv C \big( G \big)$. Hence, $\textbf{P}_h \big[ A \big( v \big) \big] \leq \mathrm{exp} \big( - c \big| v \big| \big)$. Moreover, there exist a.s. finitely many $v$ such that $A \big( v \big)$ occurs, and hence finitely many dual circuits in $G^{*}$ for which $A \big( v \big)$ occurs. The desired statement holds. \boxed{}

\bigskip

\noindent We conclude with the proof of the main theorem with the arguments below. We will argue, by contradiction with \textbf{Proposition} \textit{1}, to arrive to the conclusion that the two parameters $h_c$ and $\widetilde{h_c}$ must be equal.

\bigskip

\noindent \textit{Proof of Theorem 1}. Recall the definition of the two height parameters $h_c$ and $\widetilde{h_c}$ provided in \textit{1.3} on Page $3$. To show that $h_c = \widetilde{h_c}$, first observe that $h_c \geq \widetilde{h_c}$, because of the fact that a larger height parameter of the GFF must be taken in order for $\textbf{P}_h \big( x \overset{\geq h}{\longleftrightarrow} + \infty \big)$ to occur with positive probability. To demonstrate that the reverse inequality holds, we argue by contradiction. If $h_c \geq \widetilde{h_c}$ were to hold instead of $\widetilde{h_c} > h_c$, then there would exist intermediate height parameters, with $\widetilde{h_c} \leq \widetilde{h^1_c} \leq \widetilde{h^2_c} \leq h_c$, for which $\textbf{P}_{\widetilde{h^1_c}} \big[ \mathscr{H}\mathscr{C} \big( n , 2n \big) \big] > 0$ uniformly in $n$, by \textbf{Corollary} \textit{1}. However, because $\widetilde{h^1_c} \leq h_c$ by assumption, $\textbf{P}_{\widetilde{h_c}} \big[  0 \overset{\geq h}{\longleftrightarrow} \partial \Lambda_n \big] \leq \textbf{P}_{\widetilde{h^1_c}}\big[  0 \overset{\geq h}{\longleftrightarrow} \partial \Lambda_n  \big] \leq \textbf{P}_{\widetilde{h^2_c}}\big[  0 \overset{\geq h}{\longleftrightarrow} \partial \Lambda_n  \big] \leq \textbf{P}_{h_c}\big[  0 \overset{\geq h}{\longleftrightarrow} \partial \Lambda_n  \big] $. As $n \longrightarrow + \infty$, $\textbf{P}_{h_c}\big[  0 \overset{\geq h}{\longleftrightarrow} \partial \Lambda_n  \big] \longrightarrow 0$, in which case $\textbf{P}_{h_1} \big[  \mathscr{H}\mathscr{C} \big( n , 2n \big) \big] \longrightarrow 1$ as $n \longrightarrow + \infty$, by \textbf{Corollary} \textit{2}.

\bigskip

\noindent Over the dual graph, there exists a \textit{dual configuration} $\omega^{*}$ for which $\textbf{P}_{h_1} \big[ \omega^{*} : \omega^{*} \in \mathscr{V}\mathscr{C} \big( n , 2n \big) \big] $. Applying \textbf{Proposition} \textit{1} to the \textit{dual configuration} implies,

\begin{align*}
   \textbf{P}_{h_2} \big[ u \underset{*}{\overset{\geq h}{\longleftrightarrow}} v  \big] \leq \mathrm{exp} \big( - c \big| u - v \big| \big) \text{ } \text{ , } 
\end{align*}

\noindent for any vertices $u,v \in G^{*}$. However, the fact that the inequality above holds for vertices on the dual graph contradicts a previous result, as, 

\begin{align*}
\textbf{P}_h \big[   u \underset{*}{\overset{\geq h}{\longleftrightarrow}} v   \big] \leq \mathrm{exp} \big( - c \big| u - v \big| \big)   \Longleftrightarrow       h_2 < h_c   \text{ } \text{ . } 
\end{align*}

\noindent Hence, $\widetilde{h_c} > h_c$, and $\widetilde{h_c} = h_c$, from which we conclude the argument. \boxed{}

\section{References}

\noindent [1] Ding, J. $\&$ Wirth, M. Percolation for level-sets of the Gaussian free field on metric graphs. arXiv:
1807.11117v2 (2019).

\bigskip

\noindent [2]  Duminil-Copin, H., Karrila, A., Manolescu, I., Oulamara, M. Delocalization of the height function of
the six-vertex model. \textit{arXiv: 2012.13750 v2} (2022).

\bigskip

\noindent [3] Duminil-Copin, H., Goswami, S., Raoufi, A., Severo, F., Yadin, A. Existence of phase transition for
percolation using the Gaussian free field. \textit{Duke Math. J.} \textbf{169}(18): 3539-3563 (2020).

\bigskip

\noindent [4] Duminil-Copin, H., Manolescu, I. The phase transitions of the planar random-cluster and Potts models with q $\ge$ 1 are sharp. \textit{Probab. Theory Relat. Fields} \textbf{164}: 865–892 (2016). https://doi.org/10.1007/s00440-015-0621-0.

\bigskip

\noindent [5] Duminil-Copin, H., Hongler, C., Nolin, P. Connection probabilities and RSW-type bounds for the
FK Ising Model. \textit{Communications on Pure and Applied Mathematics} \textbf{64}(9) (2011).

\bigskip

\noindent [6] Duminil-Copin, H., Sidoravicius, V., Tassion, V. Continuity of the phase transition for planar
random-cluster and Potts models for $1 \leq q \leq 4$. \textit{Communications in Mathematical Physics} \textbf{349}: 47-107
(2017).

\bigskip

\noindent [7] Duminil-Copin, H. Tassion, V. Renormalization of crossing probabilities in the planar random-cluster
model. \textit{Moscow Mathematical Journal} \textbf{20}(4):711-740 (2020).

\bigskip

\noindent [8] Grimmett, G. The random-cluster model, vol. 333 of Grundlehren der Mathematischen Wissenschaften [Fundamental Principles of Mathematical Sciences], Springer-Verlag, Berlin,
2006.

\bigskip

\noindent [9] Jerison, D., Levine, L., Sheffield, S. Internal DLA and the Gaussian free field. Duke Math. J. \textbf{163}(2): 267-308 (2011).

\bigskip

\noindent [10] Rigas, P. Renormalization of crossing probabilities in the dilute Potts model. \textit{arXiv: 2111.10979} (2022).

\bigskip

\noindent [11] Rigas, P. From logarithmic delocalization of the six-vertex height function under sloped boundary conditions to weakened crossing probability estimates for the Ashkin-Teller, generalized random-cluster, and $(q_{\sigma},q_{\tau})$-cubic models. \textit{arXiv: 2211.14934} (2022).

\bigskip

\noindent [12] Rigas, P. Kesten's incipient infinite cluster for the three-dimensional, metric-graph Gaussian free field, from critical level-set percolation, and for the Villain model, from random cluster geometries and a Swendsen-Wang type algorithm. \textit{arXiv: 2212.07749} (2022).

\bigskip

\noindent [13] Rodriguez, P-F. A 0–1 law for the massive Gaussian free field. \textit{Probab. Theory Relat. Fields} \textbf{169}:901–930 (2017). https://doi.org/10.1007/s00440-016-0743-z.

\end{document}